\documentclass[hidelinks,12pt]{amsart}
\usepackage{bbm}
\usepackage{libertine}
\usepackage{enumitem}
\usepackage[T1]{fontenc}
\usepackage{stmaryrd}
\usepackage{color}
\DeclareFontFamily{OT1}{pzc}{}
\DeclareFontShape{OT1}{pzc}{m}{it}{<-> s * [1.100] pzcmi7t}{}
\DeclareMathAlphabet{\mathpzc}{OT1}{pzc}{m}{it}\DeclareMathAlphabet{\mathcal}{OMS}{cmsy}{m}{n}
\usepackage{fullpage}
\usepackage{mdframed}
\usepackage{graphicx}
\usepackage{calligra}
\usepackage{wasysym}
\usepackage{csquotes}
\usepackage{bm}
\usepackage{amsmath}
\usepackage[all]{xy}
\usepackage{amssymb}
\usepackage{amsthm}
\usepackage{mathrsfs}
\usepackage[OT2,T1]{fontenc}
\usepackage{centernot}
\usepackage{hyperref}

\newcommand{\DMO}[2]{\DeclareMathOperator{#1}{#2}}

\DMO{\BK}{BK}
\DMO{\FL}{FL}
\DMO{\Ann}{Ann}
\DMO{\std}{std}
\DMO{\antidiag}{antidiag}
\DMO{\locadm}{loc.adm}
\DMO{\Inj}{Inj}
\DMO{\LL}{LL}
\DMO{\Dmod}{\emph{D }-mod}
\DMO{\univ}{univ}
\DMO{\Fitt}{Fitt}
\DMO{\WD}{WD}
\DMO{\geom}{geom}
\DMO{\Fl}{Fl}
\DMO{\grad}{grad}
\DMO{\labmda}{\lambda}
\DMO{\Iw}{Iw}
\DMO{\tor}{tor}
\DMO{\coh}{coh}
\DMO{\vol}{vol}
\DMO{\semsim}{ss}
\DMO{\free}{free}
\DMO{\Alg}{Alg}
\DMO{\oth}{otherwise}
\DMO{\Ber}{Ber}
\DMO{\Diff}{Diff}
\DMO{\br}{br}
\DMO{\Isot}{Isot}
\DMO{\prim}{prim}
\DMO{\RAH}{RAH}
\DMO{\Sets}{Sets}
\DMO{\cone}{cone}
\DMO{\Grps}{Grps}
\DMO{\Dec}{Dec}
\DMO{\Flat}{Flat}
\DMO{\AbGps}{AbGps}
\DMO{\Sch}{Sch}
\DMO{\AH}{AH}
\DMO{\cl}{cl}
\DMO{\sk}{sk}
\DMO{\HC}{HC}
\DMO{\cosk}{sk}
\DMO{\ur}{ur}
\DMO{\LocSys}{LocSys}
\DMO{\rk}{rk}
\DMO{\NT}{NT}
\DMO{\cork}{cork}
\DMO{\KS}{KS}
\DMO{\MU}{MU}
\DMO{\der}{der}
\DMO{\Art}{Art}
\DMO{\Proj}{Proj}
\DMO{\End}{End}
\DMO{\Betti}{Betti}
\DMO{\Sym}{Sym}
\DMO{\cInd}{cInd}
\DMO{\GL}{GL}
\DMO{\Gal}{Gal}
\DMO{\Br}{Br}
\DMO{\Der}{Der}
\DMO{\Sp}{Sp}
\DMO{\Tan}{Tan}
\DMO{\Spin}{Spin}
\DMO{\Var}{Var}
\DMO{\Nrd}{Nrd}
\DMO{\cusp}{cusp}
\DMO{\Mat}{Mat}
\DMO{\Isom}{Isom}
\DMO{\Stab}{Stab}
\DMO{\SO}{SO}
\DMO{\Res}{Res}
\DMO{\Lie}{Lie}
\DMO{\SU}{SU}
\DMO{\Ad}{Ad}
\DMO{\ad}{ad}
\DMO{\im}{im}
\DMO{\Frob}{Frob}
\DMO{\Fr}{Fr}
\DMO{\red}{red}
\DMO{\an}{an}
\DMO{\Pic}{Pic}
\DMO{\Tor}{Tor}
\DMO{\Hdg}{Hdg}
\DMO{\id}{id}
\DMO{\pr}{pr}
\DMO{\Mor}{Mor}
\DMO{\Ext}{Ext}
\DMO{\ML}{ML}
\DMO{\PGL}{PGL}
\DMO{\SL}{SL}
\DMO{\GU}{GU}
\DMO{\GSp}{GSp}
\DMO{\GSL}{GSL}
\DMO{\Aff}{Aff}
\DMO{\NS}{NS}
\DMO{\gr}{gr}
\DMO{\Ch}{Ch}
\DMO{\QCoh}{QCoh}
\DMO{\Coh}{Coh}
\DMO{\inv}{inv}
\DMO{\Gr}{Gr}
\DMO{\Bun}{Bun}
\DMO{\Hk}{Hk}
\DMO{\GH}{GH}
\DMO{\HT}{HT}
\DMO{\LT}{LT}
\DMO{\Int}{Int}
\DMO{\UU}{U}
\DMO{\OO}{O}
\DMO{\Loc}{Loc}
\DMO{\Conn}{Conn}
\DMO{\sing}{sing}
\DMO{\si}{si}
\DMO{\Sen}{Sen}
\DMO{\MaxSpec}{MaxSpec}
\DMO{\ran}{ran}
\DMO{\coker}{coker}
\DMO{\DIV}{div}
\DMO{\Cl}{Cl}
\DMO{\Frac}{Frac}
\DMO{\VEC}{Vec}
\DMO{\Weil}{Weil}
\DMO{\SPLIT}{split}
\DMO{\Tr}{Tr}
\DMO{\val}{val}
\DMO{\pv}{p.v.}
\DMO{\disc}{disc}
\DMO{\trdeg}{tr.deg}
\DMO{\rad}{rad}
\DMO{\codim}{codim}
\DMO{\dist}{dist}
\DMO{\length}{length}
\DMO{\diam}{diam}
\DMO{\Supp}{Supp}
\DMO{\Ass}{Ass}
\DMO{\ord}{ord}
\DMO{\RE}{Re}
\DMO{\Sh}{Sh}
\DMO{\IM}{Im}
\DMO{\Tot}{Tot}
\DMO{\Bl}{Bl}
\DMO{\lcm}{lcm}
\DMO{\ann}{ann}
\DMO{\arcsinh}{arcsinh}
\DMO{\CHAR}{char}
\DMO{\MOD}{mod}
\DMO{\BB}{BB}
\DMO{\new}{new}
\DMO{\alg}{alg}
\DMO{\Irr}{Irr}
\DMO{\res}{res}
\DMO{\rank}{rank}
\DMO{\naive}{naive}
\DMO{\tors}{tors}
\DMO{\Perf}{Perf}
\DMO{\Sht}{Sht}
\DMO{\Perv}{Perv}
\DMO{\soc}{soc}
\DMO{\Mod}{Mod}
\DMO{\cyc}{cyc}
\DMO{\SC}{sc}
\DMO{\SP}{sp}
\DMO{\Deck}{Deck}
\DMO{\PSL}{PSL}
\DMO{\Area}{Area}
\DMO{\Cont}{Cont}
\DMO{\sgn}{sgn}
\DMO{\Cat}{Cat}
\DMO{\Cov}{Cov}
\DMO{\rig}{rig}
\DMO{\FSch}{FSch}
\DMO{\Rig}{Rig}
\DMO{\Spv}{Spv}
\DMO{\Spa}{Spa}
\DMO{\trace}{trace}
\DMO{\cont}{cont}
\DMO{\aff}{aff}
\DMO{\cor}{cor}
\DMO{\CH}{CH}
\DMO{\Spec}{Spec}
\DMO{\rec}{rec}
\DMO{\LGC}{LGC}
\DMO{\un}{un}
\DMO{\conj}{conj}
\DMO{\Eval}{Eval}
\DMO{\JH}{JH}
\DMO{\can}{can}
\DMO{\Fss}{Fss}
\DMO{\Speh}{Speh}
\DMO{\Ind}{Ind}
\DMO{\ch}{ch}
\DMO{\nr}{nr}
\DMO{\Swan}{Swan}
\DMO{\St}{St}
\DMO{\Ho}{Ho}
\DMO{\HH}{HH}
\DMO{\trop}{trop}
\DMO{\Jac}{Jac}
\DMO{\vir}{vir}
\DMO{\coll}{coll}
\DMO{\reg}{reg}
\DMO{\dlog}{dlog}
\DMO{\Div}{Div}
\DMO{\ab}{ab}
\DMO{\Tam}{Tam}
\DMO{\Ran}{Ran}
\DMO{\IC}{IC}
\DMO{\Sat}{Sat}
\DMO{\Rat}{Rat}
\DMO{\loc}{loc}
\DMO{\ev}{ev}
\DMO{\st}{st}
\DMO{\pst}{pst}
\DMO{\Fil}{Fil}
\DMO{\cris}{cris}
\DMO{\dR}{dR}
\DMO{\Rep}{Rep}
\DMO{\Sel}{Sel}
\DMO{\spec}{spec}
\DMO{\Spf}{Spf}
\DMO{\JL}{JL}
\DMO{\BGL}{BGL}
\DMO{\Arc}{Arc}
\DMO{\MHS}{MHS}
\DMO{\Nm}{Nm}
\DMO{\holim}{holim}
\DMO{\nInd}{nInd}
\DMO{\sSets}{s\textbf{Sets}}
\DMO{\sArt}{s\textbf{Art}}
\DMO{\BDJ}{BDJ}
\DMO{\GV}{GV}
\DMO{\BM}{BM}
\DMO{\Ord}{Ord}
\DMO{\mult}{mult}
\DMO{\WDRep}{WDRep}
\DMO{\Aut}{Aut}
\DMO{\Hom}{Hom}
\DMO{\sph}{sph}
\DMO{\Def}{Def}
\DMO{\GO}{GO}
\DMO{\diag}{diag}
\DMO{\cond}{cond}
\DMO{\ind}{ind}
\DMO{\irr}{irr}
\DMO{\RHom}{RHom}
\DMO{\sm}{sm}
\DMO{\sss}{ss}
\DMO{\sHom}{sHom}
\DMO{\Tran}{Tran}
\DMO{\Rees}{Rees}
\DMO{\lcv}{lcv} 
\DMO{\SN}{SN}
\DMO{\triv}{triv}
\DMO{\height}{ht}
\DMO{\proj}{proj}
\DMO{\Fun}{Fun}
\DMO{\cts}{cts}
\DMO{\Obj}{Obj}
\DMO{\Sing}{Sing}
\DMO{\Pro}{Pro}
\DMO{\Ig}{Ig}
\DMO{\Ha}{Ha}
\DMO{\BC}{BC}
\DMO{\RZ}{RZ}
\DMO{\supp}{supp}
\DMO{\projdim}{proj.dim}
\DMO{\Zar}{Zar}
\DMO{\Ban}{Ban}
\DMO{\LA}{LA}
\DMO{\ess}{ess}
\DMO{\op}{op}
\DMO{\Func}{Func}
\DMO{\Born}{Born}
\DMO{\Comm}{Comm}
\DMO{\Dr}{Dr}
\DMO{\LC}{LC}
\DMO{\nind}{n-ind}
\DMO{\perf}{perf}
\DMO{\charpoly}{char.poly}

\makeatletter
\def\thmhead@plain#1#2#3{%
  \thmname{#1}\thmnumber{\@ifnotempty{#1}{ }\@upn{#2}}%
  \thmnote{ {\the\thm@notefont#3}}}
\let\thmhead\thmhead@plain
\makeatother

\newtheorem*{thm1*}{Theorem}
\newtheorem*{lemma*}{Lemma}
\newtheorem*{defn1*}{Definition}

\newtheorem{thm2}{Theorem}[section]
\newtheorem{lem2}[thm2]{Lemma}

\newtheorem*{prop*}{Proposition}

\newtheorem{prop2}[thm2]{Proposition}

\newtheorem*{conj*}{Conjecture}

\theoremstyle{definition}

\newtheorem{defn2}[thm2]{Definition}

\newtheorem*{defn2*}{Definition}

\newtheorem{homework}{}
\newtheorem*{prb*}{Problem}

\newtheorem*{claim*}{Claim}

\newtheorem{rmk2}[thm2]{Remark}

\newtheorem*{rmk*}{Remark}

\newtheorem{exam2}[thm2]{Example}

\newtheoremstyle{theoremdd}
  {6pt}
  {6pt}
  {}
  {0pt}
  {\bfseries}
  {.}
  { }
  {\thmname{#1}\thmnumber{ #2}\textnormal{\thmnote{ #3}}}
  \theoremstyle{theoremdd}

\newtheoremstyle{theoremee}
  {6pt}
  {6pt}
  {\itshape}
  {0pt}
  {\bfseries}
  {.}
  { }
  {\thmname{#1}\thmnumber{ #2}\textnormal{\thmnote{ #3}}}
  \theoremstyle{theoremee}

\newcommand{\xrar}[1]{\xrightarrow{#1}}

\newcommand{\riso}{\xrar{\sim}}

 \newenvironment{psmat}
  {\left(\begin{smallmatrix}}
  {\end{smallmatrix}\right)}

 \newenvironment{psmatrix}
  {\left(\begin{smallmatrix}}
  {\end{smallmatrix}\right)}
  
 \newenvironment{pmat}
  {\begin{pmatrix}}
  {\end{pmatrix}}

\newcommand{\undertext}[2]{\mathrel{\underset{\makebox{\text{\normalfont\tiny\sffamily#1}}}{#2}}}

\newcommand{\wt}{\widetilde}
\newcommand{\wh}{\widehat}

\newcommand{\ov}{\overline}

\newcommand{\rar}{\rightarrow}

\newcommand{\ncom}[1]{\newcommand{#1}}
\ncom{\sbuset}{\subset}

\newcommand{\hrar}{\hookrightarrow}
\newcommand{\thrar}{\twoheadrightarrow}

\newcommand{\emphC}[1]{\textsf{\textbf{#1}}}

\DeclareSymbolFont{cyrletters}{OT2}{wncyr}{m}{n}
\DeclareMathSymbol{\Sha}{\mathalpha}{cyrletters}{"58}
\DeclareMathSymbol{\CyrE}{\mathalpha}{cyrletters}{"03}

\newcommand{\bs}{\backslash}
\newcommand{\et}{\operatorname{\acute{e}t}}

\DMO{\bmr}{\mathbbm{r}}
\DMO{\bmf}{\mathbbm{f}}
\DMO{\bmx}{\mathbbm{x}}

\newcommand{\bA}{\mathbb{A}}
\newcommand{\bB}{\mathbb{B}}
\newcommand{\bC}{\mathbb{C}}

\newcommand{\bN}{\mathbb{N}}

\newcommand{\bP}{\mathbb{P}}
\newcommand{\bQ}{\mathbb{Q}}

\newcommand{\bZ}{\mathbb{Z}}

\newcommand{\sH}{\mathscr{H}}

\newcommand{\cC}{\mathcal{C}}

\newcommand{\cH}{\mathcal{H}}

\newcommand{\cL}{\mathcal{L}}
\newcommand{\cM}{\mathcal{M}}

\newcommand{\cO}{\mathcal{O}}

\newcommand{\fm}{\mathfrak{m}}

\newcommand{\MSRI}{\let\thefootnote\relax\footnotetext{A part of the material is based upon work supported by the National Science Foundation under Grant No. DMS-1928930 while the author was in residence at the Mathematical Sciences Research Institute in Berkeley, California during the Spring 2023 semester.}}

\DMO{\str}{str}
\DMO{\rel}{rel}
\DMO{\Crys}{Crys}
\DMO{\Map}{Map}
\DMO{\Homeo}{Homeo}
\DMO{\AMod}{AMod}
\DMO{\Dol}{Dol}
\DMO{\sub}{sub}
\DMO{\quo}{quo}
\DMO{\Vect}{Vect}
\DMO{\ac}{ac}
\DMO{\PGSp}{PGSp}
\DMO{\Ssp}{sp}

\DMO{\Out}{Out}
\DMO{\Mer}{Mer}
\DMO{\comp}{comp}
\DMO{\Inv}{Inv}
\DMO{\MIC}{MIC}
\DMO{\Heis}{Heis}
\DMO{\Isoc}{Isoc}
\DMO{\la}{la}
\DMO{\corank}{corank}
\DMO{\Td}{Td}
\DMO{\hol}{hol}
\DMO{\gen}{gen}
\DMO{\Conf}{Conf}
\DMO{\OConf}{OConf}
\DMO{\metab}{metab}

\begin{document}
\title{Arithmetic quantum local systems over the moduli of curves}
\author{Gyujin Oh}\address{Department of Mathematics, Columbia University, 2990 Broadway, New York, NY 10027}
\maketitle
\begin{abstract}
We construct an arithmetic analogue of the quantum local systems on the moduli of curves, and study its basic structure. Such an arithmetic local system gives rise to a uniform way of assigning a Galois cohomology class of the first geometric \'etale cohomology of a smooth proper curve over a number field.
\end{abstract}
\tableofcontents
\section{Introduction}
The \emph{quantum representations} of mapping class groups $\Mod_{g,n}$ are widely used terms referring to the (projective) representations arising from ($3$-dimensional) topological quantum field theories (TQFTs) or those arising from the (projectively) flat connections on the vector bundle of conformal blocks (the \emph{Hitchin connections}). They provide an abundance of representations of the mapping class groups that do not factorize through the Torelli groups, and are useful for the representation theory of the mapping class groups. For example, certain families of quantum representations are known to be \emph{asymptotically linear}, i.e., the intersection of the kernel of all representations in a family is trivial (\cite{Andersen}, \cite{FWW}). Many quantum representations are also known to be integral (\cite{Gilmer}), and they can be used to study the modular representation theory of mapping class groups and symplectic groups (\cite{GM}).

In this paper, we construct the arithmetic analogues of certain quantum representations of mapping class groups. Recall that the mapping class group $\Mod_{g,n}$ is the topological fundamental group of the moduli stack $\cM_{g,n,\bC}$ of compact Riemann surfaces of genus $g$ with $n$ marked points (see \cite{Oda}). Thus, a representation of $\Mod_{g,n}$ may be considered as a topological local system on $\cM_{g,n,\bC}$. In this paper, we will construct the analogue of the quantum representations constructed in \cite{BPS} (see also \cite{DeRenziMartel}) as $p$-adic \'etale local systems on $\cM_{g,n,\bQ}$. We hope that one could construct the arithmetic analogue of more general quantum local systems  using similar techniques, if their homological construction is known to exist (cf. \cite{BFS}).

More precisely, what we achieve in this paper is as follows. Let $K$ be a number field, and let $Q\in\cM_{g,1}(K)$ be a $K$-rational point, corresponding to a smooth proper curve $C$ of genus $g$ over $K$ (with a $K$-rational point). Then, for an odd prime $p$, we construct the \emph{artihmetic Heisenberg local system} $\rho_{Q,p}^{\Heis}$, which is an \'etale $\bZ_{p}$-local system of rank $2g+1$ over $\cM:=\cM_{g,1,K}$. This is Definition \ref{def:arithHeis}. The construction crucially relies on the relative version of the Puiseux section construction of  Anderson, Ihara, and Matsutomo (e.g., \cite{Matsumoto}). From this, one may construct the arithmetic analogues of quantum local systems by composing with a Weil representation.

The arithmetic Heisenberg local system is a non-split extension of the $\bZ_{p}$-linear dual of the relative $H^{1}_{\et}$ of the universal curve $\cC\rar\cM$ by the cyclotomic character,
\[0\rar \bZ_{p}(1)\rar\rho_{Q,p}^{\Heis}\rar \sH^{1}_{\et}(\cC/\cM,\bZ_{p})^{\vee}\rar0.\]
This is shown in Lemma \ref{lem:nonsplit} and Proposition \ref{prop:subquo}. Restricting the \'etale local system at a $K$-rational point $x\in\cM(K)$, corresponding to a smooth proper curve $C'$ of genus $g$ over $K$ (with a $K$-rational point), we obtain a Galois cohomology class
\[c_{Q,x}^{\Heis}\in H^{1}_{\cL^{S\cup\lbrace p\rbrace}}(K,H^{1}_{\et}(C'_{\ov{K}},\bZ_{p})(1)),\]where $S$ is the set of primes of $K$ at which either $C$ or $C'$ has bad reduction, and $\cL^{S\cup\lbrace p\rbrace}$ is the local condition of being unramified at each prime of $K$ away from $S\cup\lbrace p\rbrace$ (see Proposition \ref{prop:HeisenbergExtension})\footnote{We also expect that $c_{Q,x}^{\Heis}$ is crystalline when both $C$ and $C'$ have good reduction at $p$; see Remark \ref{rem:crystalline}.}. Thus, the arithmetic Heisenberg local system $\rho_{Q,p}^{\Heis}$ may be regarded as a systematic way of obtaining a class in the Selmer group of the first geometric \'etale cohomology of a smooth proper curve of genus $g$ over $K$, which may have arithmetic applications.

The extension class $c_{Q,x}^{\Heis}$ depends crucially on both $Q$ and $x$ (cf. \cite{Mochizuki}, \cite{Faltings}). Depending on $Q$ and $x$, $c_{Q,x}^{\Heis}$ may be simply zero. For example, we suspect that $c_{Q,x}^{\Heis}=0$ when $Q=x$, but $c_{Q,x}^{\Heis}$ when $Q\ne x$ seems more mysterious, due to the nature of the relative Puiseux section construction. We expect that $c_{Q,x}^{\Heis}$ and its generalizations are arithmetically interesting. For example, in the case of genus $0$, the conformal block local system is related to the Knizhnik--Zamolodchikov (KZ) equations, and their homological/motivic realizations were constructed in \cite{Looijenga} and \cite{BBM}. It was recently shown in \cite{BFM} that the arithmetic analogue of the KZ local system, restricted at certain rational points, yields interesting Galois representations (referred to as \emph{KZ motives}). We hope to be able to come back to the problem of computing $c_{Q,x}^{\Heis}$ in future work.

%
%
\subsection{Acknowledgements}
We thank Richard Wentworth and Junho Peter Whang for inspiring suggestions.
\subsection{Notations}
Given a group $G$, $\wh{G}$ is the profinite completion of $G$, and $\wh{G}_{p}$ is the pro-$p$ completion of $G$. A profinite group is finitely generated if there is a dense subgroup generated by finitely many elements. 

A compact surface of genus $g$ is denoted as $\Sigma_{g}$, and $\Sigma_{g,n}$ is $\Sigma_{g}$ with $n$ points removed. 
\section{Relative Puiseux sections of homotopy exact sequences}\label{RelativePuiseux}
Our goal of this preliminary section is to develop a relative variant of the theory of tangential basepoints, or \emph{Puiseux sections}, as developed by Deligne, Anderson, Ihara, and Matsumoto (\cite{Deligne}, \cite{AndersonIhara}, \cite{Ihara}, \cite{Matsumoto}) to produce a particular splitting of the homotopy exact sequence of \'etale fundamental groups. Below, we recall the relevant works in \cite{Matsumoto} and extend them to the relative setting.

Our setup is as follows. Let $K\subset\bC$ be a subfield of $\bC$, and let $\ov{K}\subset\bC$ be the algebraic closure of $K$ in $\bC$. Let $f:X\rar S$ be a dominant morphism, such that $S$ is smooth over $K$, and $X$ is an integral $K$-variety. Let $D\subset X$ be a relative normal crossings divisor over $S$ (see \cite[XIII.2.1]{SGA1} for the definition), and let $Y=X-D$, which is smooth over $S$. Let $P\in D(K)$, and let $Q\in S(K)$ be the image of $P$. Suppose that  $X\rar S$ is smooth at $P$. 
Choose a small enough simply connected subset $U\subset S(\bC)$ of $Q$ (which may  not be open nor closed). Let $\cM_{S}^{\an}(U)$ be the ring of germs of meromorphic functions in open neighborhoods around $U$ which can be analytically continued into finitely multivalued unramified meromorphic functions on the whole $S(\bC)$. Let $\cM_{S}(U)\subset\cM_{S}^{\an}(U)$ be the algebraic closure of the field of rational functions $K(S)$ in $\cM_{S}^{\an}(U)$.

Let $\pi_{1}(S(\bC),U)$ be the topological fundamental group whose basepoint is any point\footnote{The fundamental group is canonically identified upon any choice of the basepoint in $U$, as $U$ is simply connected.} of $U$. Then, by \cite[Proposition 1.2]{Matsumoto}, the group homomorphism $\pi_{1}(S(\bC),U)\rar \Aut\cM_{S}(U)$, where a loop $\gamma$ maps to the ``analytic continuation along $\gamma$'' map $\gamma^{\ac}:\cM_{S}(U)\rar\cM_{S}(U)$, gives rise to  an isomorphism $\pi_{1}(S(\bC),U)^{\wedge}\riso\Gal(\cM_{S}(U)/\ov{K}(S))$, and $\cM_{S}(U)$ is a maximal algebraic extension of $\ov{K}(S)$ unramified over $S\otimes_{K}\ov{K}$. Let $\eta_{S}:\Spec K(S)\rar S$ be the generic point of $S$, and let $\ov{\eta_{S}}:\Spec\ov{K(X)}\rar S$ be a geometric point underlying $\eta_{S}$, where $K(S)$ is regarded as a subfield of $K(X)$. By \cite[V.8.2]{SGA1}, the natural homomorphism $\Gal(\ov{K(S)}/K(S))=\pi_{1,\et}(\eta_{S},\ov{\eta_{S}})\rar\pi_{1,\et}(S,\ov{\eta_{S}})$ is surjective and factors through an isomorphism $\Gal(\cM_{S}(U)/K(S))\riso\pi_{1,\et}(S,\ov{\eta_{S}})$.

Let $\eta:\Spec K(X)\rar X$ be the generic point of $X$, and let $\ov{\eta}:\Spec\ov{K(X)}\rar X$ be a geometric point underlying $\eta$. Note that these are also the generic point and the geometric generic point of $Y$ as well. Also, note that $\ov{\eta_{S}}=f\circ\ov{\eta}$. Our goal is to construct a splitting of the homotopy exact sequence
\begin{equation}\label{HES}\tag{$*$}1\rar\pi_{1,\et}(Y_{\ov{\eta_{S}}},\ov{\eta})\rar\pi_{1,\et}(Y,\ov{\eta})\rar\pi_{1,\et}(S,\ov{\eta_{S}})\rar1,\end{equation}using the relative version of the Puiseux section as in \cite[Proposition 1.3]{Matsumoto}.

\begin{prop2}\label{prop:RelativePuiseux}
Let $d=\dim S$ and $e=\dim X$. Let $z_{1},\cdots,z_{d}\in\fm_{S,Q}\subset\cO_{S,Q}$ form a local coordinate system of $S$ at $Q$, and let $w_{1},\cdots,w_{e}\in\fm_{X,P}\subset\cO_{X,P}$ form a  local coordinate system of $X$ at $P$, such that $w_{i}=f^{*}z_{i}$ for $1\le i\le d$. Suppose that $D$ is locally around $P$ cut out by the equation $w_{c}\cdots w_{e}=0$ for some $d<c<e$. Suppose also that the closed subscheme $S'\subset X$ cut out by $w_{d+1}=\cdots=w_{e}=0$ satisfies the property that, on the  reduced subscheme $(S'_{0})_{\red}$ of the irreducible component $S'_{0}\subset S'$ containing $P$,  $f\rvert_{(S'_{0})_{\red}}:(S'_{0})_{\red}\riso S$ is an isomorphism.

Let $U$ be a  small open neighborhood of $P$ in $X(\bC)$ such that $w_{1},\cdots,w_{e}$ exist as holomorphic functions on $U$, with $w_{i}(P)=0$, $1\le i\le e$, and $(w_{1},\cdots,w_{e}):U\riso D(0,\varepsilon)^{e}$ is a bijection onto the $e$-th power of a one-dimensional complex open disc of radius $\varepsilon>0$ for some small $\varepsilon$. We define
\[V:=\lbrace x\in U:w_{i}(x)\in(0,\varepsilon)\text{ for $c\le i\le e$}, w_{i}(x)=0\text{ for other $i$}\rbrace,\]so that $(w_{1},\cdots,w_{e}):V\riso\lbrace0\rbrace^{c-1}\times(0,\varepsilon)^{e-c+1}$ is bijective. On a small neighborhood of $V$, for any $N\in\bN$ and $c\le i\le e$, we may define a univalued meromorphic function $w_{i}^{1/N}$ by 
\[w_{i}^{1/N}(x):=e^{\frac{\log w_{i}(x)}{N}}.\]
Let $\cM_{Y}^{f,P}(V)$ be the integral closure of $\cO_{Y,P}\subset K(Y)=K(X)$ in $\cM_{X}(V)$.
\begin{enumerate}
\item Each $h\in\cM_{Y}^{f,P}(V)$ has a convergent relative Puiseux expansion
\[h=\sum_{i_{d+1},\cdots,i_{e}\in\bN}a_{i_{d+1},\cdots,i_{e}}w_{d+1}^{i_{d+1}}\cdots w_{c-1}^{i_{c-1}}w_{c}^{i_{c}/N}\cdots w_{e}^{i_{e}/N},\quad a_{i_{d+1},\cdots,i_{e}}\in\cM_{S}(Q),\]for some $N\in\bN$.
\item If we let $\sigma\in\Gal(\cM_{S}(Q)/K(S))$ act on the coefficients of $h$ to obtain $\sigma(h)$, then $\sigma(h)\in\cM_{Y}^{f,P}(V)$. Thus, this induces a splitting
\[s_{P}:\pi_{1,\et}(S,\ov{\eta_{S}})\cong\Gal(\cM_{S}(Q)/K(S))\rar\Gal(\cM_{Y}(V)/K(Y))\cong\pi_{1,\et}(Y,\ov{\eta}),\]of the homotopy exact sequence \emph{(}\ref{HES}\emph{)}.
\end{enumerate}
All of the above remain true when $X$ and $S$ are  \emph{Deligne--Mumford stacks}\footnote{We use the definition that an algebraic stack is \emph{connected} (\emph{irreducible}, respectively) if its underlying topological space is connected (irreducible, respectively).}  over $K$.
\end{prop2}
\begin{proof}
The following argument is an adaptation of the proof of \cite[Proposition 1.3]{Matsumoto} in a relative setting.

Note that $w_{1},\cdots,w_{e}$ are rational functions on $X$. As $\pi_{1}(U-D(\bC))\cong\bZ^{c-e+1}$, there exists an $N\in\bN$ such that, on the cover defined by the Cartesian squares
\[\xymatrix{Y'\ar@{^(->}[rr]\ar[d] && X'\ar[rr]\ar[d] && (\bP^{1})^{c-e+1}\ar[d]^-{(z_{1},\cdots,z_{c-e+1})\mapsto(z_{1}^{N},\cdots,z_{c-e+1}^{N})} 
\\ Y\ar@{^(->}[rr] && X\ar[rr]_-{(w_{c},\cdots,w_{e})}&&(\bP^{1})^{c-e+1}.}\]
Let $U'=U\times_{X(\bC)}X'(\bC)$, and let $W$ be an open neighborhood of $Q$ in $S(\bC)$ such that $(z_{1},\cdots,z_{d}):W\riso D(0,\varepsilon)^{d}$ is a bijection. Note that there is a unique point $P'\in U'$ that sits over $P\in U$. Then, $h$, as a finitely multivalued unramified meromorphic function on $Y'(\bC)\cap U'$, has no monodromy around the divisors $w_{i}^{1/N}=0$ for $c\le i\le e$, so $h$ extends to a univalued meromorphic function on $U'$, regular at $P'$.  By the definition of the local coordinates, $X'\rar X\xrar{f} S$ restricts to $U'\rar U\rar W$, which we again denote as $f$ by abuse of notation. On $U'$, the functions $w_{1},\cdots,w_{c-1},w_{c}^{1/N},\cdots,w_{e}^{1/N}$ exist as univalued holomorphic functions, and \[(f,w_{d+1},\cdots,w_{c-1},w_{c}^{1/N},\cdots,w_{e}^{1/N}):U'\riso W\times D(0,\varepsilon)^{e-d},\] is a bjiection. Note that, on $X'$, $w_{1},\cdots,w_{c-1},w_{c}^{1/N},\cdots,w_{e}^{1/N}$ are rational functions. Taking the Taylor series expansion of $h$ on $U'$ around the origin gives rise to an expansion of the form
\[h=\sum_{i_{d+1},\cdots,i_{e}\in\bN}a_{i_{d+1},\cdots,i_{e}}w_{d+1}^{i_{d+1}}\cdots w_{c-1}^{i_{c-1}}w_{c}^{i_{c}/N}\cdots w_{e}^{i_{e}/N},\]where each $a_{i_{d+1},\cdots,i_{e}}$ is a meromorphic function on $W$ holomorphic at $Q$. By \cite[Lemma 1.1]{Matsumoto} ($K$ in the lemma is $K(S)$ in our context), $a_{i_{d+1},\cdots,i_{e}}\in\ov{K(S)}$. 

We claim that each $a_{i_{d+1},\cdots,i_{e}}$ is realized as a rational function on a finite \'etale cover of $S$, thereby implying that $a_{i_{d+1},\cdots,i_{e}}\in\cM_{S}(Q)$. As $h$ is a rational function on a finite \'etale cover of $Y$, we may regard it as a rational function on a finite \'etale cover of $Y'$, as $Y'\rar Y$ is a finite \'etale cover. Note that $\bP^{1}\xrar{z\mapsto z^{N}}\bP^{1}$ restricts to $\lbrace0\rbrace\riso \lbrace0\rbrace$, so the closed subscheme $S''\subset X'$ cut out by $w_{d+1}=\cdots=w_{c-1}=w_{c}^{1/N}=\cdots=w_{e}^{1/N}=0$ is isomorphic to $S'$. Let $S''_{0}$ be the irreducible component of $S''$ containing $P'$. Then, $f\rvert_{(S''_{0})_{\red}}:(S''_{0})_{\red}\riso S$ is an isomorphism. As $h$ has no monodromy around $w_{i}^{1/N}=0$, $c\le i\le e$, it is unramified over $S''$. As $a_{0,\cdots,0}=h\rvert_{S''}$ and $h$ is unramified over $S''$, $a_{0,\cdots,0}$ is a rational function on a finite \'etale cover of $S''$, thus on a finite \'etale cover of $S$. 

Suppose that $h$ is a rational function on $X''$, where $X''\rar X'$ is a finite  cover, unramified over $Y''=X'-D'$, where $D=\left(\bigcup_{c\le i\le e}\lbrace w_{i}^{1/N}=0\rbrace\right)\cup D'$. As $w_{d+1},\cdots,w_{c-1},w_{c}^{1/N},\cdots,w_{e}^{1/N}$ are rational functions on $X''$, we see that $\frac{\partial h}{\partial w_{d+1}},\cdots,\frac{\partial h}{\partial w_{c-1}},\frac{\partial h}{\partial w_{c}^{1/N}},\cdots,\frac{\partial h}{\partial w_{e}^{1/N}}$ are again rational functions on $X''$ holomorphic at the preimage of $P$. Thus, each such partial derivative is an element of $\cM_{Y}^{f,P}(V)$. Therefore, by induction on $i_{d+1}+\cdots+i_{e}$, we can conclude that $a_{i_{d+1},\cdots,i_{e}}\in\cM_{S}(Q)$. This proves (1).

Before we prove (2), we first note that, inside $\ov{K(X)}$, $\ov{K(S)}\cap K(X)=K(S)$. Indeed, if $x\in\ov{K(S)}\cap K(X)$, then $f(x)=0$ for some monic irreducible polynomial $f(T)\in K(S)[T]$. On the other hand, due to the section $(f\rvert_{S'})^{-1}:S\riso (S'_{0})_{\red}\hrar X$, we have a field morphism $p:K(X)\rar K(S)$ such that $p(y)=y$ for $y\in K(S)\subset K(X)$. Therefore, $f(x)=0$ implies $f(p(x))=0$, so $f(T)=T-p(x)$, which means $x=p(x)\in K(S)$.

From the proof of (1), we know that $a_{i_{d+1},\cdots,i_{e}}\in K(S''')$, where $S'''=X''\times_{X'}S''$. We may assume that $K(S''')/K(S)$ is Galois.  Our goal is to show that, for any $\sigma\in\Gal(\cM_{S}(Q)/K(S))$, $\sigma(h)\in\cM_{Y}^{f,P}(V)$. As $\cM_{S}(Q)/K(S''')/K(S)$ is a Galois subextension, $\sigma$ acts on the coefficients $a_{i_{d+1},\cdots,i_{e}}$ via the image of $\sigma$ along the quotient $\Gal(\cM_{S}(Q)/K(S))\thrar\Gal(K(S''')/K(S))$. We may just let the image be denoted $\sigma\in\Gal(K(S''')/K(S))$ for simplicity.

Let $L=K(S''')\otimes_{K(S)}K(X)$, which is a field as $K(S''')\cap K(X)=K(S)$ inside $\ov{K(X)}$. Let $f(T)\in L[T]$ be the monic minimal polynomial of $h$ over $L$. Let $M=L[T]/(f^{\sigma}(T))$, where $f^{\sigma}(T)\in L[T]$ is the polynomial obtained by applying $\sigma\otimes\id:K(S''')\otimes_{K(S)}K(X)\rar K(S''')\otimes_{K(S)}K(X)$ on the coefficients. Let $M'$ be the Galois closure of $M$, and let $\varphi:\wt{X}\rar X$ be the normalization of $X$ in $M'$. Then, $\varphi$ factors through $\wt{X}\xrar{\varphi_{1}} X\times_{S}S'''\xrar{\varphi_{2}} X$, $T$ is a univalued function on $\wt{X}$, and both $\varphi^{-1}(V)\subset \wt{X}(\bC)$ and $\varphi_{2}^{-1}(V)\subset X\times_{S}S'''(\bC)$ are disjoint unions of  (analytically) isomorphic copies of $V$. Let $V'\subset\varphi_{2}^{-1}(V)$ be a connected component of $\varphi_{2}^{-1}(V)$. Restricting $T$ to each connected component of $\varphi_{1}^{-1}(V')$, we get a function in $\cM_{Y\times_{S}S'''}^{f,P'}(V')\subset\cM_{Y}^{f,P}(V)$, where $P'\in V'$ is sent to $P\in V$ via $\varphi_{2}:V'\riso V$. Note that there are $[M':L]$ many connected components of $\varphi_{1}^{-1}(V')$. If $T$ restricts to the same function in $\cM_{Y\times_{S}S'''}^{f,P'}(V')$ on two different connected components of $\varphi_{1}^{-1}(V')$, then there exists a nontrivial element $\tau\in\Gal(M'/L)$ corresponding to a path between two components on which $T$ restricts to the same function, which implies that $\tau(T)=T$. Thus, $\tau$ fixes $M$. Therefore, we get $\frac{|\Gal(M'/L)|}{|\Gal(M'/M)|}=[M:L]$ different functions in $\cM_{Y\times_{S}S'''}^{f,P'}(V')$ by restricting $T$ to different components of $\varphi_{1}^{-1}(V')$. As $[M:L]=\deg f$, these functions are exactly the solutions of $f^{\sigma}(T)$ in the formal Puiseux series with coefficients in $\cM_{S}(Q)$. As $\sigma(h)$ is one such solution, $\sigma(h)\in\cM_{Y\times_{S}S'''}^{f,P'}(V')\subset\cM_{Y}^{f,P}(V)$. This gives rise to a homomorphism $s_{P}:\Gal(\cM_{S}(Q)/K(S))\rar\Gal(\cM_{Y}(V)/K(Y))$. Now note that the surjective map $\pi_{1,\et}(Y,\ov{\eta})\rar\pi_{1,\et}(S,\ov{\eta_{S}})$ is identified with $\Gal(\cM_{Y}(V)/K(Y))\thrar\Gal(\cM_{S}(Q)\otimes_{K(S)}K(Y)/K(Y))$ (here $\cM_{S}(Q)\otimes_{K(S)}K(Y)$ is a subfield of $\cM_{S}(Q)$), and the composition of this with $s_{P}$ is obviously an isomorphism by the definition of $s_{P}$. Therefore, $s_{P}$ provides a splitting to the homotopy exact sequence (\ref{HES}).

The above argument applies verbatim in the case when $X$ and $S$ are  Deligne--Mumford stacks, by defining holomorphic functions on orbifolds  to be holomorphic functions on an \'etale chart satisfying compatibility under the corresponding \'etale equivalence relation and using the orbifold fundamental groups.
\end{proof}

\section{Arithmetic Heisenberg local systems}
Let $g\ge2$ and $n\ge1$ be integers. Let $\cM_{g,n}$ be the moduli of curves of genus $g$ with $n$  marked points over $\bQ$, which is a smooth geometrically irreducible Deligne--Mumford stack over $\bQ$. Namely, $\cM_{g,n}\rar\Sch_{\bQ,\et}$ is the category fibered in groupoids over the category of $\bQ$-schemes with \'etale topology, such that, over a $\bQ$-scheme $S$, $(\cM_{g,n})_{S}$ is the groupoid of $S$-families of genus $g$ curves with $n$ marked points, in the following sense.
\begin{defn2}An \emph{$S$-family of genus $g$ curves with $n$ marked points} is a tuple $(C,\sigma_{1},\cdots,\sigma_{n})$, where $f: C\rar S$ is a smooth proper morphism and $\sigma_{1},\cdots,\sigma_{n}$ are sections of $f$ such that, for every point $s:\Spec k\rar S$, $C_{s}$ is a genus $g$ connected curve and $\sigma_{1}(s),\cdots,\sigma_{n}(s)\in C_{s}$ are distinct points.
\end{defn2}It is well-known that, for any $S$-family of genus $g$ curves with $n$ marked points $(C,\sigma_{1},\cdots,\sigma_{n})$, the image of $\sigma_{i}$ is a divisor of $C$. For such a family, we define $C^{\circ}\subset C$ to be the open subscheme obtained by removing the divisors $\cup_{i=1}^{n}\sigma_{i}(S)$. We define the \emph{ordered configuration space} $\OConf_{m}(C^{\circ}/S)$ of $m$ distinct points on $C^{\circ}$ as
\[\OConf_{m}(C^{\circ}/S):=(\undertext{$m$ times}{\underbrace{C^{\circ}\times_{S}C^{\circ}\times_{S}\cdots\times_{S}C^{\circ}}})\bs\Delta,\]where $\Delta$ is the diagonal divisor (i.e. the union of all divisors with repeated entries). Note that $S_{m}$ acts on $\OConf_{m}(C^{\circ}/S)$ as the permutation of the entries. We then define the (unordered) \emph{configuration space} $\Conf_{m}(C^{\circ}/S)$ of $m$ distinct points on $C^{\circ}$ as
\[\Conf_{m}(C^{\circ}/S):=\OConf_{m}(C^{\circ}/S)/S_{m},\] where the he quotient is taken in the category of fppf sheaves over $S$, which is an $S$-scheme by \cite[Tag 07S7]{stacks}, as $\OConf_{m}(C^{\circ}/S)\rar S$ is quasi-affine and the action of $S_{m}$ is free.
\begin{lem2}\label{Geometry}
\hfill
\begin{enumerate}
\item The formation of $\Conf_{m}$ and $\OConf_{m}$ commutes with the base-change.
\item For any $S$-family of genus $g$ curves with $n$ marked points, the morphism $\Conf_{m}(C^{\circ}/S)\rar S$ is smooth.
\item For any $S$-family of genus $g$ curves with $n$ marked points, if $S$ is smooth over $\bQ$ and irreducible, $\Conf_{m}(C^{\circ}/S)$ is irreducible.
\end{enumerate}
\end{lem2}
\begin{proof}
(1) is immediate from the definition of the configuration space as an fppf quotient. For (2), note that $\OConf_{m}(C^{\circ}/S)\rar\Conf_{m}(C^{\circ}/S)$ is an $S_{m}$-torsor (again by the proof of \cite[Tag 07S7]{stacks}), so is in particular \'etale. As smoothness is \'etale-local on the source, and as $C^{\circ}\rar S$ is smooth, $\Conf_{m}(C^{\circ}/S)\rar S$ is smooth. For (3), it suffices to show that $C^{\circ}\times_{S}\cdots\times_{S}C^{\circ}$ is irreducible. As $C^{\circ}\times_{S}\cdots\times_{S}C^{\circ}$ is $\bQ$-smooth, it suffices to show that it is connected. Suppose not, so that it is expressed as a union of two disjoint nonempty open subsets $U_{1},U_{2}$. As $\pi:C^{\circ}\times_{S}\cdots\times_{S}C^{\circ}\rar S$ is smooth, it is open, so $\pi(U_{1}),\pi(U_{2})$ are open subsets of $S$, whose union is again $S$. As $S$ is connected, $\pi(U_{1})\cap\pi(U_{2})\ne\emptyset$. Choose a point $p\in \pi(U_{1})\cap \pi(U_{2})$. Then, $\pi^{-1}(p)$ is a union of two disjoint open sets $U_{1}\cap \pi^{-1}(p)$ and $U_{2}\cap\pi^{-1}(p)$, and both are nonempty. Therefore, $\pi^{-1}(p)$ is disconnected. Thus, to show (3), it suffices to show that $C^{\circ}\times_{S}\cdots\times_{S}C^{\circ}$ is connected when $S$ is $\Spec$ of a field of characteristic zero. In that case the underlying topological space of $C^{\circ}\times_{S}\cdots\times_{S}C^{\circ}$ is a self-product of the underlying topological space of $C^{\circ}$. As $C^{\circ}$ is connected, any self-product is also connected, as desired.
\end{proof}
Let $\phi:\cC_{g,n}\rar\cM_{g,n}$ be the universal curve, and  $\Sigma_{1},\cdots,\Sigma_{n}:\cM_{g,N}\rar\cC_{g,n}$ be the universal sections. As $\Conf_{m}(\cC_{g,n}^{\circ}/\cM_{g,n})$ is irreducible by Lemma \ref{Geometry}(3), we can take its generic point $\eta_{C}$ and geometric generic point $\ov{\eta_{C}}$. Let $\eta_{M}$ and $\ov{\eta_{M}}$ be the generic point and the geometric generic point of $\cM_{g,n}$, respectively. Let $C^{\circ}_{\ov{\eta}}$ be the fiber of $\cC_{g,n}^{\circ}$ over $\ov{\eta_{M}}$. 
We have the homotopy exact sequence 
\begin{equation}\label{HESCM}\tag{$**$}1\rar\pi_{1,\et}(\Conf_{m}(C^{\circ}_{\ov{\eta}}/\ov{\eta_{M}}),\ov{\eta_{C}})\rar\pi_{1,\et}(\Conf_{m}(\cC_{g,n}^{\circ}/\cM_{g,n}),\ov{\eta_{C}})\rar\pi_{1,\et}(\cM_{g,n},\ov{\eta_{M}})\rar1.\end{equation}
\subsection{Relative Puiseux section on $\Conf_{m}(\cC^{\circ}_{g,n}/\cM_{g,n})\rar\cM_{g,n}$}\label{RelativePuiseuxMgn}
Using the relative Puiseux sections developed in \S\ref{RelativePuiseux}, we construct a splitting of the homotopy exact sequence (\ref{HESCM}). Let $K\subset\bC$ be a subfield, and $Q\in\cM_{g,n}(K)$ be a $K$-rational point, corresponding to $(C,\sigma_{1},\cdots,\sigma_{n})$. Let $z_{1},\cdots,z_{3g+n-3}\in\fm_{\cM_{g,n},Q}\subset\cO_{\cM_{g,n},Q}$ form a local coordinate system of $\cM_{g,n}$ at $Q$. Let $w_{0},w_{1},\cdots,w_{3g+n-3}\in\fm_{\cC_{g,n},\Sigma_{1}(Q)}$ form a local coordinate system of $\cC_{g,n}$ at $\Sigma_{1}(Q)$, where $w_{i}=\phi^{*}z_{i}$ for $1\le i\le 3g+n-3$, and the divisor $\Sigma_{1}(\cM_{g,n})$ is locally cut by the equation $w_{0}=0$. Note that $w_{0}$ is sent to a uniformizer of $\fm_{C,\sigma_{1}}$ via the natural quotient map. We define $\ov{\cC_{g,n}^{\circ}}\rar\cM_{g,n}$ as $\ov{\cC_{g,n}^{\circ}}:=\cC_{g,n}-\bigcup_{i=2}^{n}\Sigma_{i}(\cM_{g,n})$.

Using $w_{0}$, we obtain a morphism $w:\cC_{g,n,K}\rar\bP^{1}_{K}$, which is smooth of dimension $3g+n-3$ at $\Sigma_{1}(Q)$. Note that $(w,\phi):\cC_{g,n,K}\rar\bP^{1}_{K}\times_{K}\cM_{g,n,K}=\bP^{1}_{\cM_{g,n,K}}$ is \'etale at $\Sigma_{1}(Q)$. As in \cite[\S1.4]{Matsumoto}, we define a morphism
\[\alpha:\bA^{m}_{t}:=\Spec K[t_{1},\cdots,t_{m}]\rar\bA^{m}_{u}:=\Spec K[u_{1},\cdots,u_{m}],\]by $u_{i}\mapsto t_{i}t_{i+1}\cdots t_{m}$. 
Then $\alpha$ restricts to an isomorphism
\[\alpha:\bA^{m}_{t}-D_{t}\riso \OConf_{m}(\bA^{1}_{K}-\lbrace0\rbrace/\Spec K),\]
where $D_{t}$ is the union of divisors $\lbrace t_{i}=0\rbrace$, for $1\le i\le m-1$, $\lbrace t_{i}=1\rbrace$, for $1\le i\le m-1$, $\lbrace t_{m}=0\rbrace$, and $\lbrace t_{i}t_{i+1}\cdots t_{j}=1\rbrace$, for $1\le i<j\le m-1$. We denote $w^{m}$ the restriction of $w\times w\times\cdots \times w$ to $\cC_{m}$, where \[\cC_{m}:=\undertext{$m$ times}{\underbrace{\left(\ov{\cC^{\circ}_{g,n,K}}-w^{-1}(\infty)\right)\times_{\cM_{g,n,K}}\cdots\times_{\cM_{g,n,K}}\left(\ov{\cC^{\circ}_{g,n,K}}-w^{-1}(\infty)\right)}}.\]Consider the Cartesian square
\[\xymatrix{X\ar[r]^-{w_{X}}\ar[d]&\bA^{m}_{t}\ar[d]^-{\alpha}\\\cC_{m}\ar[r]_-{w^{m}}&\bA^{m}_{u}.}\]
Note that $(w^{m},\phi):\cC_{m}\rar\bA^{m}_{u}\times_{K}\cM_{g,n,K}$ is \'etale at $(\Sigma_{1}(Q),\cdots,\Sigma_{1}(Q))$. Therefore, $(w_{X},\phi):X\rar\bA^{m}_{t}\times_{K}\cM_{g,n,K}$ is \'etale at the fiber in $X$ above $(\Sigma_{1}(Q),\cdots,\Sigma_{1}(Q))$, and this fiber maps isomorphically via $(w_{X},\phi)$ to the fiber in $\bA^{n}_{t}\times_{K}\cM_{g,n,K}$ over $(0,\cdots,0,Q)\in\bA^{m}_{u}\times_{K}\cM_{g,n,K}$. Thus, there is a unique point $P$ on this fiber which maps to $(0,\cdots,0,Q)\in\bA^{m}_{u}\times_{K}\cM_{g,n,K}$ and $P$ is $K$-rational. We define $Y$ by the fiber product
\[\xymatrix{Y\ar@{^(->}[r]
\ar[d]&
\cC_{m}\ar[d]^-{w^{m}}\\ 
\OConf_{m}(\bA^{1}_{K}-\lbrace0\rbrace/\Spec K)\ar@{^(->}[r]
&\bA_{u}^{m}\\
\bA_{t}^{m}-D_{t}\ar[u]^-{\cong}\ar@{^(->}[r]
&\bA_{t}^{m}\ar[u]_-{\alpha}.}\]
Clearly $Y$ is an open subspace of $\OConf_{m}(\cC_{g,n}^{\circ}/\cM_{g,n})$. Moreover, we have the following Cartesian diagram,
\[\xymatrix{Y\ar@{^(->}[r]\ar[d]
&X\ar[d]^-{w_{X}}\\\bA^{m}_{t}-D_{t}\ar@{^(->}[r]
&\bA_{t}^{m},}\]which implies that the complement $D$ of $Y$ in $X$ is a normal crossing divisor locally at $P$, cut out locally by $t_{1}\cdots t_{m}=0$, where $t_{1},\cdots,t_{m}$ are the local coordinates  of $X$ at $P$ pulled back from that of $\bA^{m}_{t}$ at $(0,\cdots,0)$. 
\begin{lem2}\label{lem:GeometryX}\hfill
\begin{enumerate}
\item The morphism $X\rar \cM_{g,n,K}$ is dominant, and $X$ is integral.
\item The smooth locus of the morphism $X\rar \cM_{g,n,K}$ contains $Y$ and $P$.
\item The closed subscheme $S'\subset X$ cut out by $t_{1}=\cdots=t_{m}=0$ satisfies the property that $X\rar\cM_{g,n,K}$ restricts to $(S'_{0})_{\red}\riso\cM_{g,n,K}$, where $S'_{0}$ is the reduced subscheme of the irreducible component of $S'$ containing $P$. 
\end{enumerate}
\end{lem2}
\begin{proof}
Note that an open dense subspace of $X$ is isomorphic to an open dense subspace of \[\undertext{$m$ times}{\underbrace{\ov{\cC^{\circ}_{g,n,K}}\times_{\cM_{g,n,K}}\cdots\times_{\cM_{g,n,K}}\ov{\cC^{\circ}_{g,n,K}}}}.\]As the map from this space to $\cM_{g,n,K}$ is a smooth surjective morphism with irreducible fibers, and as $\cM_{g,n,K}$ is irreducbile, it follows that the above space is irreducible (e.g. \cite[Tag 004Z]{stacks}). It follows that $X\rar \cM_{g,n,K}$ is dominant, and that $X$ is irreducible. 

Let $\bigcup_{j\in J}V_{j}$ be an \'etale chart of $\cM_{g,n,K}$ consisting of affine schemes, and over each $V_{j}$, $\bigcup_{i\in I_{j}}U_{ij}$ be an affine open cover of $\ov{\cC^{\circ}_{g,n,K}}-w^{-1}(\infty)\times_{\cM_{g,n,K}}V_{j}$. Let $V_{j}=\Spec B_{j}$ and $U_{ij}=\Spec A_{ij}$. Restricting $w$ to $U_{ij}$ corresponds to an element $f_{ij}\in A_{ij}$. Then, $X$ is covered by the affine schemes of the form
\[\Spec \frac{\left(A_{i_{1}j}\otimes_{B_{j}}\cdots\otimes_{B_{j}}A_{i_{m}j}\right)[u_{1},\cdots,u_{m-1}]}{((f_{i_{1}}\otimes\cdots\otimes1)u_{1}-(1\otimes f_{i_{2}}\otimes\cdots\otimes1),\cdots,(1\otimes\cdots\otimes f_{i_{m-1}}\otimes 1)u_{m-1}-(1\otimes\cdots\otimes f_{i_{m}}))}\]
\[=\Spec\left(A_{i_{1}j}\otimes_{B_{j}}\cdots\otimes_{B_{j}}A_{i_{m}j}\right)\left\langle\frac{1\otimes f_{i_{2}}\otimes\cdots\otimes1}{f_{i_{1}}\otimes\cdots\otimes1},\cdots,\frac{1\otimes\cdots\otimes f_{i_{m}}}{1\otimes\cdots\otimes f_{i_{m-1}}\otimes 1}\right\rangle,\]for $i_{1},\cdots,i_{m}\in I$, where $\left(A_{i_{1}j}\otimes_{B_{j}}\cdots\otimes_{B_{j}}A_{i_{m}j}\right)\left\langle\frac{1\otimes f_{i_{2}}\otimes\cdots\otimes1}{f_{i_{1}}\otimes\cdots\otimes1},\cdots,\frac{1\otimes\cdots\otimes f_{i_{m}}}{1\otimes\cdots\otimes f_{i_{m-1}}\otimes 1}\right\rangle$ is the subring of $\Frac\left(A_{i_{1}j}\otimes_{B_{j}}\cdots\otimes_{B_{j}}A_{i_{m}j}\right)$ generated over $\left(A_{i_{1}j}\otimes_{B_{j}}\cdots\otimes_{B_{j}}A_{i_{m}j}\right)$ by the specified elements. This ring is evidently reduced, so $X$ is integral. This proves (1).

As $Y$ is an open subspace of $\cC_{m}$, $Y\rar\cM_{g,n,K}$ is smooth. We have already seen that $X\rar\cM_{g,n,K}$ is smooth at $P$. This proves (2).

Note that $(0,\cdots,0)\in\bA_{t}^{m}$ is sent isomorphically to $(0,\cdots,0)\in\bA_{u}^{m}$ via $\alpha$, so $S'$ is sent isomorphically to $\lbrace w=0\rbrace^{m}$ via $X\rar\cC_{m}$. Note that, in a Zariski open neighborhood of $P$, $\lbrace w=0\rbrace^{m}$ coincides with $\Sigma_{1}(\cM_{g,n})\times_{\cM_{g,n}}\cdots\times_{\cM_{g,n}}\Sigma_{1}(\cM_{g,n})\cong\cM_{g,n}$. As $\cM_{g,n}$ is irreducible, $S'_{0}$ and $\cM_{g,n}$ underlie the same topological space. Therefore, $(S'_{0})_{\red}\rar\cM_{g,n}$ is an isomorphism. This proves (3).
\end{proof}
Using Lemma \ref{lem:GeometryX} and that $(w_{X},\phi):X\rar\bA^{n}_{t}\times_{K}\cM_{g,n,K}$ is \'etale at $P$, we can apply Proposition \ref{prop:RelativePuiseux} to $Y\rar\cM_{g,n,K}$, with $P$ and $Q$. Thus, we obtain the relative Puiseux section
\[s_{P}:\pi_{1,\et}(\cM_{g,n,K},\ov{\eta_{M}})\rar\pi_{1,\et}(Y,\ov{\eta_{OC}}),\]of the natural surjective morphism $\pi_{1,\et}(Y,\ov{\eta_{OC}})\rar\pi_{1,\et}(\cM_{g,n,K},\ov{\eta_{M}})$, where $\eta_{OC}$ and $\ov{\eta_{OC}}$ are the generic point and the geometric generic point of $\OConf_{m}(\cC_{g,n,K}^{\circ}/\cM_{g,n,K})$, respectively. As the natural surjective morphism $\pi_{1,\et}(Y,\ov{\eta_{OC}})\rar\pi_{1,\et}(\cM_{g,n,K},\ov{\eta_{M}})$ factors through $\pi_{1,\et}(Y,\ov{\eta_{OC}})\rar\pi_{1,\et}(\OConf_{m}(\cC_{g,n,K}^{\circ}/\cM_{g,n,K}),\ov{\eta_{OC}})\rar\pi_{1,\et}(\Conf_{m}(\cC_{g,n,K}^{\circ}/\cM_{g,n,K}),\ov{\eta_{C}})\rar\pi_{1,\et}(\cM_{g,n,K},\ov{\eta_{M}})$, the relative Puiseux section $s_{P}$ also gives rise to the relative Puiseux section 
\[s_{P}:\pi_{1,\et}(\cM_{g,n,K},\ov{\eta_{M}})\rar\pi_{1,\et}(\Conf_{m}(\cC_{g,n,K}^{\circ}/\cM_{g,n,K}),\ov{\eta_{C}}),\]of the natural surjective morphism $\pi_{1,\et}(\Conf_{m}(\cC_{g,n,K}^{\circ}/\cM_{g,n,K}),\ov{\eta_{C}})\rar\pi_{1,\et}(\cM_{g,n,K},\ov{\eta_{M}})$, by the abuse of notation. 


\subsection{Arithmetic Heisenberg local systems}
Using the relative Puiseux section $s_{P}$ constructed in \S\ref{RelativePuiseuxMgn}, we define the arithmetic analogue of the Heisenberg local systems defined in \cite{BPS}. 
From now on, we always assume that $n=1$, as was done in \emph{op.~cit}. 
The action of $\pi_{1,\et}(\cM_{g,1,K},\ov{\eta_{M}})$ on $\pi_{1,\et}(\Conf_{m}(C_{\ov{\eta}}^{\circ}/\ov{\eta_{M}}),\ov{\eta_{C}})$ via the section $s_{P}$ gives rise to a homomorphism
\[\rho_{Q}:\pi_{1,\et}(\cM_{g,1,K},\ov{\eta_{M}})\rar\Aut_{\cont}\pi_{1,\et}(\Conf_{m}(C_{\ov{\eta}}^{\circ}/\ov{\eta_{M}}),\ov{\eta_{C}}),\]where $\Aut_{\cont}$ is the group of continuous automorphisms. 
By \cite[Corollary XII.5.2]{SGA1} and \cite[Theorem 1.1]{Landesman}, $\pi_{1,\et}(\Conf_{m}(C_{\ov{\eta}}^{\circ}/\ov{\eta_{M}}),\ov{\eta_{C}})\cong \pi_{1,}(\Conf_{m}(C^{\circ}(\bC)/\bC),*)^{\wedge}$ for any embedding of $\bQ(\cM_{g,1})\hrar\bC$. In the literature, $\pi_{1}(\Conf_{m}(C^{\circ}(\bC)/\bC),*)$ is called the \emph{surface braid group}, and will be denoted $\bB_{m}(\Sigma_{g,1})$. On the other hand, $\pi_{1,\et}(\cM_{g,1,K},\ov{\eta_{M}})$ is called the \emph{arithmetic mapping class group} over $K$, and will be denoted $\AMod_{g,1,K}$ (if $K=\bQ$ we will often omit $K$ from the subscript). The name originates from the fact that $\AMod_{g,1,K}$ sits in between the short exact sequence
\[1\rar\pi_{1,\et}(\cM_{g,1,\ov{K}},\ov{\eta_{M}})\rar\AMod_{g,1,K}\rar\Gal(\ov{K}/K)\rar1,\]and the first term of the short exact sequence is identified with the profinite completion of the mapping class group $\Mod_{g,1}^{\wedge}$ by applying the same comparison results between \'etale and topological fundamental groups over algebraically closed fields of characteristic $0$ to \cite{Oda}, which shows that  $\pi_{1}(\cM_{g,1}(\bC),*)$ is isomorphic to the mapping class group $\Mod_{g,1}$.


Thus, from the relative Puiseux section $s_{P}$, we obtain a homomorphism
\[\rho_{Q}:\AMod_{g,1,K}\rar\Aut_{\cont}\left(\bB_{m}(\Sigma_{g,1})^{\wedge}\right).\]
This is in general a very complicated homomorphism;  $\Aut_{\cont}(\bB_{m}(\Sigma_{g,1})^{\wedge})$ is a massive group. For example, it is known that $\Out_{\cont}\bB_{m}(\Sigma_{0,1})^{\wedge}$ contains the Grothendieck--Teichm\"uller group (see \cite{MinamideNakamura} for the detail). 

On the other hand, the metabelianization $\bB_{m}(\Sigma_{g,1})^{\metab}$ has the structure of a Heisenberg group (e.g. \cite{BGG}). Recall that, for a group $G$, the \emph{metabelianization} $G^{\metab}$ is defined by $G^{\metab}:=G/[G,[G,G]]$. It is a \emph{metabelian group}, i.e. a group whose commutator group is abelian. The metabelianization $G^{\metab}$ has the universal property that any group homomorphism $G\rar H$ to a metabelian group $H$ factors through the metabelianization $G\rar G^{\metab}\rar H$. The same definition works well for finitely generated profinite groups by \cite[Theorem 1.4]{NikolovSegal}\footnote{The correct definition of the metabelianization of a topological group $G$ would be $G^{\metab}:=G/\ov{[G,\ov{[G,G]}]}$. By \emph{loc. cit.}, we do not need to take the closures of the commutator groups.}, which satisfies the similar universal property for homomorphisms to metabelian profinite groups.

We will from now on restrict to the $p$-part of $\rho_{Q}$, where $p>2$ is an odd prime number. Namely, we consider the homomorphism
\[\rho_{Q,p}:\AMod_{g,1,K}\rar\Aut_{\cont}(\bB_{m}(\Sigma_{g,1})^{\wedge}_{p}),\]obtained by composing $\rho_{Q}$ with the natural map\footnote{One way to see the existence of the natural homomorphism is as follows. As $\bB_{m}(\Sigma_{g,1})^{\wedge}$ is a finitely generated profinite group, by \cite[Theorem 1.1]{NikolovSegal}, $\bB_{m}(\Sigma_{g,1})^{\wedge}_{p}$ is the pro-$p$ completion of $\bB_{m}(\Sigma_{g,1})^{\wedge}$, from which the existence of the natural homomorphism is obvious.} $\Aut_{\cont}(\bB_{m}(\Sigma_{g,1})^{\wedge})\rar\Aut_{\cont}(\bB_{m}(\Sigma_{g,1})^{\wedge}_{p})$. 
\begin{lem2}\label{lem:metabelian}
Any (continuous)\footnote{The continuity condition is unnecessary thanks to \cite[Theorem 1.1]{NikolovSegal}.} automorphism of $\bB_{m}(\Sigma_{g,1})^{\wedge}_{p}$ uniquely induces a (continuous) automorphism of $(\bB_{m}(\Sigma_{g,1})^{\metab})^{\wedge}_{p}$, giving a natural group homomorphism $\Aut_{\cont}(\bB_{m}(\Sigma_{g,1})^{\wedge}_{p})\rar\Aut_{\cont}((\bB_{m}(\Sigma_{g,1})^{\metab})^{\wedge}_{p})$.
\end{lem2}
\begin{proof}
For a finitely generated group $G$, we claim that $(\wh{G}_{p})^{\metab}\cong(G^{\metab})^{\wedge}_{p}$. The Lemma will then follow, as the metabelianization is clearly a (topological) characteristic quotient for finitely generated pro-$p$ groups. Note that the pro-$p$ completion of a metabelian group is metabelian (i.e. any two commutators commute with each other), so $(G^{\metab})^{\wedge}_{p}$ is a metabelian pro-$p$ group. Thus, the natural homomorphism $\wh{G}_{p}\rar (G^{\metab})^{\wedge}_{p}$ factors through $(\wh{G}_{p})^{\metab}\rar(G^{\metab})^{\wedge}_{p}$. On the other hand, the natural homomorphism $G\rar\wh{G}_{p}$ gives rise to a homomorphism $G^{\metab}\rar(\wh{G}_{p})^{\metab}$. As $(\wh{G}_{p})^{\metab}$ is pro-$p$, the universal property of the pro-$p$ completion gives a natural homomorphism $(G^{\metab})^{\wedge}_{p}\rar(\wh{G}_{p})^{\metab}$. It is straightforward to check that these two homomorphisms are inverses to each other.
%
\end{proof}
\begin{defn2}\label{def:arithHeis}
As per Lemma \ref{lem:metabelian}, we may construct the \emph{arithmetic Heisenberg local system} on the moduli of curves,
\[\rho_{Q,p}^{\Heis}:\AMod_{g,1,K}\rar\Aut_{\cont}((\bB_{m}(\Sigma_{g,1})^{\metab})_{p}^{\wedge}).\]
\end{defn2}
The reason why the above homomorphism is called the arithmetic Heisenberg local system is as follows: $(\bB_{m}(\Sigma_{g,1})^{\metab})_{p}^{\wedge}$ has a structure of a $p$-adic Heisenberg group.
\begin{lem2}\label{lem:Metab is Heisenberg}Let $m\ge3$ and $g\ge1$. 
The pro-$p$ group $(\bB_{m}(\Sigma_{g,1})^{\metab})_{p}^{\wedge}$ is the central extension of $H_{1}(\Sigma_{g},\bZ_{p})$ by $\bZ_{p}$,
\[0\rar\bZ_{p}\rar(\bB_{m}(\Sigma_{g,1})^{\metab})_{p}^{\wedge}\rar H_{1}(\Sigma_{g},\bZ_{p})\rar0,\]corresponding to the intersection symplectic pairing $\omega:H_{1}(\Sigma_{g},\bZ_{p})\times H_{1}(\Sigma_{g},\bZ_{p})\rar\bZ_{p}$.
\end{lem2}
\begin{proof}
The metabelianization $\bB_{m}(\Sigma_{g,1})^{\metab}$ is well-known (e.g. \cite[Proposition 3.13]{BGG}); it is the central extension of $H_{1}(\Sigma_{g},\bZ)$ by $\bZ$, corresponding to the intersection pairing $H_{1}(\Sigma_{g},\bZ)\times H_{1}(\Sigma_{g},\bZ)\rar\bZ$. The Lemma follows from the fact that pro-$p$ completion is exact on the category of finitely generated nilpotent groups (e.g. \cite[Exercise 1.21]{Segal}).%
%
\end{proof}
\begin{defn2}
Analogous to \cite{BPS} (see also \cite{DeRenziMartel}), we may define an \emph{arithmetic quantum local system} as follows. Note that $H_{p^{N}}:=\frac{(\bB_{m}(\Sigma_{g,1})^{\metab})_{p}^{\wedge}}{((\bB_{m}(\Sigma_{g,1})^{\metab})_{p}^{\wedge})^{p^{N}}}$ is the Heisenberg group for $H_{1}(\Sigma_{g},\bZ/p^{N}\bZ)$. By the Stone--von Neumann theorem for finite Heisenberg groups (e.g., \cite[Proposition 2.1]{Lysenko}), given a nondegenerate character $\psi:\bZ/p^{N}\bZ\rar\ov{\bQ}_{p}$ (i.e, $\psi(x)=\zeta_{p^{N}}^{x}$ for a primitive $p^{N}$-th root of unity $\zeta_{p^{N}}\in\ov{\bQ}_{p}$), there is a unique irreducible representation $(\rho_{\psi},V_{\psi})$ of $H_{p^N}$ of central character $\psi$ over $\ov{\bQ}_{p}$. This gives rise to a projective representation $\Aut(H_{p^{N}})\rar\PGL(V_{\psi})$, which can be composed with $\rho_{Q,p}^{\Heis}$ to define a representation $\AMod_{g,1,K}\rar\PGL(V_{\psi})$. 
\end{defn2}
The (continuous) automorphism group of the $p$-adic Heisenberg group $(\bB_{m}(\Sigma_{g,1})^{\metab})_{p}^{\wedge}$ is manageable.
\begin{lem2}\label{lem:AutH}The continuous automorphism group $\Aut_{\cont}((\bB_{m}(\Sigma_{g,1})^{\metab})_{p}^{\wedge})$ can be identified with a parabolic subgroup of $\GL_{2g+1}(\bZ_{p})=\GL(\bZ_{p}\oplus H_{1}(\Sigma_{g},\bZ_{p}))$,
\[\Aut_{\cont}((\bB_{m}(\Sigma_{g,1})^{\metab})_{p}^{\wedge})\cong\begin{pmat}\GL_{1}(\bZ_{p})&* \\ 0&\GSp_{2g}(\bZ_{p})\end{pmat}\subset\GL_{2g+1}(\bZ_{p}).\]
\end{lem2}
\begin{proof}
Let $\Aut_{\cont}((\bB_{m}(\Sigma_{g,1})^{\metab})_{p}^{\wedge})^{0}$ be the subgroup of $\Aut_{\cont}((\bB_{m}(\Sigma_{g,1})^{\metab})_{p}^{\wedge})$ consisting of the elements that fix the commutator subgroup $\bZ_{p}$ of $(\bB_{m}(\Sigma_{g,1})^{\metab})_{p}^{\wedge}$. Then, 
\[\frac{\Aut_{\cont}((\bB_{m}(\Sigma_{g,1})^{\metab})_{p}^{\wedge})}{\Aut_{\cont}((\bB_{m}(\Sigma_{g,1})^{\metab})_{p}^{\wedge})^{0}}\cong\GL_{1}(\bZ_{p}).\]
Thus, it suffices to show that
\[\Aut_{\cont}((\bB_{m}(\Sigma_{g,1})^{\metab})_{p}^{\wedge})^{0}\cong\begin{pmat}1&* \\ 0&\GSp_{2g}(\bZ_{p})\end{pmat}\subset\GL_{2g+1}(\bZ_{p}).\]Here $\GSp_{2g}(\bZ_{p})$ is really regarded as $\GSp(H_{1}(\Sigma_{g},\bZ_{p}),\omega)$.
Note that the unipotent radical is naturally identified with $\Hom_{\bZ_{p}}(H_{1}(\Sigma_{g},\bZ_{p}),\bZ_{p})$. Given $\psi\in\Hom_{\bZ_{p}}(H_{1}(\Sigma_{g},\bZ_{p}),\bZ_{p})$ and $g\in\GSp(H_{1}(\Sigma_{g},\bZ_{p}))$, we define $\varphi_{\psi,g}\in\Aut_{\cont}((\bB_{m}(\Sigma_{g,1})^{\metab})_{p}^{\wedge})$ as
\[\varphi_{\psi,g}(a,x):=(a+\psi(x),gx).\]It is straightforward to check that $\varphi_{\psi,g}\in\Aut_{\cont}((\bB_{m}(\Sigma_{g,1})^{\metab})_{p}^{\wedge})^{0}$ and that this gives rise to the desired isomorphism.\end{proof}
\begin{rmk2}
Note that Lemma \ref{lem:AutH} \emph{does not} mean that an element of the parabolic subgroup as defined in the Lemma acts $\bZ_{p}$-linearly. After all, $(\bB_{m}(\Sigma_{g,1})^{\metab})^{\wedge}_{p}$ is a non-abelian group.
\end{rmk2}
Taking the Levi quotient, we obtain maps \[\Aut_{\cont}((\bB_{m}(\Sigma_{g,1})^{\metab})^{\wedge}_{p})\rar\GL_{1}(\bZ_{p}),\quad \Aut_{\cont}((\bB_{m}(\Sigma_{g,1})^{\metab})^{\wedge}_{p})\rar\GSp_{2g}(\bZ_{p})\rar\GL_{2g}(\bZ_{p}).\]
\begin{defn2}
We define $\rho_{Q,\sub,p}^{\Heis}$ and $\rho_{Q,\quo,p}^{\Heis}$ be the subrepresentation and the quotient representation of $\rho_{Q,p}^{\Heis}$ corresponding to composing $\rho_{Q,p}^{\Heis}$ with the natural maps  $\Aut_{\cont}((\bB_{m}(\Sigma_{g,1})^{\metab})^{\wedge}_{p})\rar\GL_{1}(\bZ_{p})$ and $\Aut_{\cont}((\bB_{m}(\Sigma_{g,1})^{\metab})^{\wedge}_{p})\rar\GL_{2g}(\bZ_{p})$, respectively.
\end{defn2}
Thus, the (a priori non-linear) information from the metabelian group  gives us  a linear (!) extension of $\bZ_{p}$-\'etale local systems on $\cM_{g,1,K}$,
\begin{equation}
\label{eq:extension}\tag{$\dagger$}1\rar\rho_{Q,\sub,p}^{\Heis}\rar\rho_{Q,p}^{\Heis}\rar\rho_{Q,\quo,p}^{\Heis}\rar1.
\end{equation}
Note that $\rho_{Q,\sub,p}^{\Heis}$ is a character and  $\rho_{Q,\quo,p}^{\Heis}$ is of dimension $2g$. Note also that the extension is non-split, as its geometric part $\rho_{Q,p}^{\Heis}\rvert_{\Mod_{g,1}^{\wedge}}$ is indecomposable.
\begin{lem2}\label{lem:nonsplit}The extension \emph{(\ref{eq:extension})} restricted to the geometric mapping class group $\Mod_{g,1}^{\wedge}$ is non-split.
\end{lem2}
\begin{proof}
We use the notation of \cite[Fig. 1]{BG}. Note that the image of $[\delta_{1},\delta_{2}]\in\bB_{m}(\Sigma_{g,1})$ in $\bB_{m}(\Sigma_{g,1})^{\metab}$ is a nonzero element in the commutator subgroup. One may thus accordinly apply Dehn twists to $\delta_{1}$ to obtain an element in the commutator subgroup of $\bB_{m}(\Sigma_{g,1})^{\metab}$. This implies that, for any set-theoretic splitting of (\ref{eq:extension}), the $
\bZ_{p}[\Mod_{g,1}^{\wedge}]$-module generated by the lifts of the elements of $\bB_{m}(\Sigma_{g,1})^{\ab}$ in $\bB_{m}(\Sigma_{g,1})^{\metab}$ will contain $[\bB_{m}(\Sigma_{g,1})^{\metab},\bB_{m}(\Sigma_{g,1})^{\metab}]$, which implies that it would generate the whole $\rho_{Q,p}^{\Heis}\rvert_{\Mod_{g,1}^{\wedge}}$. This implies that there is no splitting respecting the $\Mod_{g,1}^{\wedge}$-action.\end{proof}
\section{Galois cohomology classes from the arithmetic Heisenberg local systems}
We first show that the local subsystem and the quotient local system are of classical nature. 
\begin{prop2}\label{prop:subquo}\hfill
\begin{enumerate}
\item The quotient local system $\rho_{Q,\quo,p}^{\Heis}$ coincides with the $p$-adic \'etale local system $(R^{1}\phi_{*,\et}\bZ_{p})^{\vee}$, where $(\cdot)^{\vee}$ is the $\bZ_{p}$-linear dual.
\item The character $\rho_{Q,\sub,p}^{\Heis}$ coincides with the composition $\AMod_{g,1,K}\thrar\Gal(\ov{K}/K)\xrar{\chi^{\cyc}_{p}}\bZ_{p}^{\times}$ where $\chi^{\cyc}_{p}$ is the $p$-adic cyclotomic character.
\end{enumerate}
\end{prop2}
\begin{proof}
We first show (1). The  quotient local system $\rho_{Q,\quo,p}^{\Heis}$ is, by definition, $(R^{1}c_{*,\et}\bZ_{p})^{\vee}$, where $c:\Conf_{m}(\cC_{g,1,K}^{\circ}/\cM_{g,1,K})\rar\cM_{g,1,K}$. 
Let $t:\Pic^{0}(\cC_{g,1,K}/\cM_{g,1,K})\rar \cM_{g,1,K}$ be the relative Picard scheme. Then, there is a natural morphism $\Conf_{m}(\cC_{g,1,K}^{\circ}/\cM_{g,1,K})\rar\Pic^{0}(\cC_{g,1,K}/\cM_{g,1,K})$, $\lbrace x_{1},\cdots,x_{m}\rbrace\mapsto x_{1}+\cdots+x_{m}-m\sigma$. As $R^{1}t_{*,\et}\bZ_{p}\cong R^{1}\phi_{*,\et}\bZ_{p}$, by functoriality, there is a natural map $R^{1}\phi_{*,\et}\bZ_{p}\cong R^{1}t_{*,\et}\bZ_{p}\rar R^{1}c_{*,\et}\bZ_{p}$, and we claim that this map is an isomorphism. As the \'etale site has enough points (see \cite[Tag 04K5]{stacks}), we may check that this is an isomorphism by checking at every geometric point of $\cM_{g,1,K}$. Therefore, it suffices to show that, for an algebraically closed field $k$ of characteristic $0$ of finite transcendence degree over $\ov{K}$ and a genus $g$ curve $C/k$ with a section $\sigma\in C(k)$, the analogous map $\Conf_{m}(C^{\circ}/k)\rar\Pic^{0}(C/k)$ induces an isomorphism $H^{1}_{\et}(\Pic^{0}(C/k),\bZ_{p})\riso H^{1}_{\et}(\Conf_{m}(C^{\circ}/k),\bZ_{p})$. We may embed $k$ into $\bC$ and the isomorphism may be checked after base-changing $k$ to $\bC$. By Artin comparison, it suffices to show that $H^{1}(\Pic^{0}(C(\bC)),\bZ_{p})\riso H^{1}(\Conf_{m}(C^{\circ}(\bC)),\bZ_{p})$, which will follow if we show that $\pi_{1}(\Conf_{m}(C^{\circ}(\bC)))^{\ab}\otimes_{\bZ}\bZ_{p}=\bB_{m}(\Sigma_{g,1})^{\ab}\otimes_{\bZ}\bZ_{p}\rar\pi_{1}(\Pic^{0}(C(\bC)))^{\ab}\otimes_{\bZ}\bZ_{p}$ is an isomorphism. By \cite[\S3.2]{BGG}, we see that the torsion group of $\bB_{m}(\Sigma_{g,1})^{\ab}$ is just $\bZ/2\bZ$, and the basis of the torsion-free part can be given by moving one of the $m$ points around the standard basis loops of $H_{1}(C(\bC),\bZ)=\pi_{1}(C(\bC))^{\ab}$. This implies that these basis elements of $\bB_{m}(\Sigma_{g,1})^{\ab}$ are sent to the images of the basis loops of $\pi_{1}(C(\bC))^{\ab}$ sent to $\pi_{1}(\Pic^{0}(C(\bC)))^{\ab}$,  translated by a particular element. This implies that $\bB_{m}(\Sigma_{g,1})^{\ab}\otimes_{\bZ}\bZ_{p}\rar\pi_{1}(\Pic^{0}(C(\bC)))^{\ab}\otimes_{\bZ}\bZ_{p}$ is an isomorphism, as desired.

We now show (2). The restriction of $\rho_{Q,p}^{\Heis}$ to $\Mod_{g,1}^{\wedge}$ comes from the relative Puiseux section of the homotopy exact sequence applied to $\Conf_{m}(\cC_{g,1,\ov{K}}^{\circ}/\cM_{g,1,\ov{K}})\rar\cM_{g,1,\ov{K}}$. Therefore, this restriction is the pro-$p$ completion of the  homomorphism $\Mod_{g,1}\rar\Aut(\pi_{1}(\Conf_{m}(C^{\circ}),*)^{\metab})$ coming from the topological fibration $\Conf_{m}(\cC_{g,1}^{\circ}(\bC)/\cM_{g,1}(\bC))\rar\cM_{g,1}$. By \cite[Proposition 13]{BPS}, it follows that $\Mod_{g,1}$ fixes the commutator subgroup $\bZ\subset\pi_{1}(\Conf_{m}(C^{\circ}),*)^{\metab}$. This implies that $\rho_{Q,\sub,p}^{\Heis}\rvert_{\Mod_{g,1}^{\wedge}}$ is trivial, so that $\rho_{Q,\sub,p}^{\Heis}$ factors through a  Galois character $\chi:\Gal(\ov{K}/K)\rar\bZ_{p}^{\times}$. 

We  use $Q$ to obtain the splitting $t_{Q}:\Gal(\ov{K}/K)\rar\pi_{1,\et}(\cM_{g,1,K},\ov{\eta_{M}})$ of the homotopy exact sequence
\[1\rar\pi_{1,\et}(\cM_{g,1,\ov{K}},\ov{\eta_{M}})\rar\pi_{1,\et}(\cM_{g,1,K},\ov{\eta_{M}})\rar\Gal(\ov{K}/K)\rar1.\]Similarly, we  use the Puiseux section associated with $P$ to obtain the splitting $t_{P}:\Gal(\ov{K}/K)\rar\pi_{1,\et}(\Conf_{m}(\cC^{\circ}_{g,1,K}/\cM_{g,1,K}),\ov{\eta_{C}})$ of the homotopy exact sequence
\[1\rar\pi_{1,\et}(\Conf_{m}(\cC^{\circ}_{g,1,\ov{K}}/\cM_{g,1,\ov{K}}),\ov{\eta_{C}})\rar\pi_{1,\et}(\Conf_{m}(\cC^{\circ}_{g,1,K}/\cM_{g,1,K}),\ov{\eta_{C}})\rar\Gal(\ov{K}/K)\rar1.\]From the construction of the relative Puiseux section $s_{P}$, it follows easily that $t_{P}=s_{P}\circ t_{Q}$. Therefore, for $g\in\Gal(\ov{K}/K)$, $\chi(g)\in\bZ_{p}^{\times}$ is such that, for any $a$ in the commutator subgroup $\bZ_{p}$ of $(\pi_{1,\et}(\Conf_{m}(C^{\circ}_{\ov{\eta}}/\ov{\eta_{M}}),\ov{\eta_{C}})^{\metab})^{\wedge}_{p}$, we have $\chi(g)a=t_{P}(g)\cdot a$. 

Note that, under the presentation used in the proof of Lemma \ref{lem:AutH}, $(b,x)(c,y)(b,x)^{-1}(c,y)^{-1}=(2\omega(x,y),0)$, regardless of what $b,c\in\bZ_{p}$ you take. Therefore, applying to $a=2\omega(x,y)$, we see that $\chi(g)$ is the scalar such that $\chi(g)\omega(x,y)=\omega(t_{Q}(g)x,t_{Q}(g)y)$, where $\pi_{1,\et}(\cM_{g,1,K},\ov{\eta_{M}})$ acts on $(\pi_{1,\et}(\Conf_{m}(C^{\circ}_{\ov{\eta}}/\ov{\eta_{M}}),\ov{\eta_{C}})^{\ab})^{\wedge}_{p}$ as the (outer) action coming from the homotopy exact sequence (\ref{HESCM}). Therefore, this is on the Galois equivariance property of $\omega$ on the restriction of the \'etale local system $(R^{1}\phi_{*,\et}\bZ_{p})^{\vee}$ at $Q$. As the intersection pairing is Poincar\'e dual to cup product pairing, and as the $G_{K}$-representation $V=H^{1}_{\et}(C_{\ov{K}},\bZ_{p})$ satisfies $V^{\vee}\cong V(1)$ for any smooth proper curve $C$ over $K$, it follows that $\omega:V^{\vee}\times V^{\vee}\rar\bZ_{p}(1)$ is Galois equivariant. This shows (2).
\end{proof}
\begin{rmk2}
One may also prove Proposition \ref{prop:subquo}(1) by applying the proof of \cite[Theorem 1.1(ii)]{Matsumoto} to the context of relative Puiseux sections.
\end{rmk2}
Therefore, given $Q\in\cM_{g,1}(K)$, the arithmetic Heisenberg local system $\rho_{Q,p}^{\Heis}$ gives a universal way of obtaining a Galois extension class of the first \'etale cohomology of a curve.
\begin{defn2}
Let $x\in\cM_{g,1}(L)$ be an $L$-rational point, for a field extension $L/K$, corresponding to a smooth proper curve $C$ over $L$, and an $L$-rational point $s\in C(L)$. By Proposition \ref{prop:subquo}, restricting (\ref{eq:extension}) at $x$ gives a short exact sequence of $\bZ_{p}[\Gal(\ov{L}/L)]$-modules,
\[0\rar \bZ_{p}(1)\rar(\rho_{Q,p}^{\Heis})\rvert_{x}\rar H^{1}_{\et}(C_{\ov{L}},\bZ_{p})(1)\rar0.\]We define the \emph{Heisenberg extension class} \[c_{Q,x}^{\Heis}\in \Ext^{1}_{\bZ_{p}[\Gal(\ov{L}/L)]}(H^{1}_{\et}(C_{\ov{L}},\bZ_{p}),\bZ_{p})= H^{1}(L,H^{1}_{\et}(C_{\ov{L}},\bZ_{p})(1)),\]as the Galois cohomology class corresponding to the (Tate twist of the) above extension. 
\end{defn2}
We end by showing that the Heisenberg extension class has the expected local properties at good primes, at least away from $p$.
\begin{prop2}\label{prop:HeisenbergExtension}
Let $L/K$ be number fields, and let $Q\in\cM_{g,1}(K)$, $x\in\cM_{g,1}(L)$. Let $C$ be the curve of genus $g$ over $K$ and $s\in C(K)$ be the rational point,  corresponding to $Q$, and let $(C',s')$ be the pair of a curve over $L$ and its $L$-rational point corresponding to $x$. 

Let $v$ be a finite prime of $L$ such that $v\nmid p$. Suppose that both $C_{L}$ and $C'$ have good reduction at $v$. Then, $c_{Q,x}^{\Heis}$ is unramified at $v$.
\end{prop2}
\begin{proof}
We may base-change $Q$ to $L$ and the problem does not change, so we may assume that $K=L$. 
We use the integral model $\cC_{g,1,\cO_{L_{v}}}\rar\cM_{g,1,\cO_{L_{v}}}$ of $\cC_{g,1,L}\rar\cM_{g,1,L}$, which is still smooth and proper. 
The construction of the relative Puiseux section carries through even for schemes over $\cO_{L_{v}}$, as long as the condition $f\rvert_{(S'_{0})_{\red}}\riso S$ being an isomorphism holds integrally. Indeed the construction of \S\ref{RelativePuiseuxMgn} works integrally. Therefore, there is a relative Puiseux section $s_{P,\cO_{L_{v}}}$ that is compatible with $s_{P}$ restricted to $L_{v}$,
\[\hspace*{-1.3cm}
\xymatrix{
1\ar[r]&
\pi_{1,\et}(\Conf_{m}(C_{\ov{\eta}}^{\circ}/\ov{\eta_{M}}),\ov{\eta_{C}})^{\wedge}_{p}\ar[r]\ar@{=}[d]&
\pi_{1,\et}'(\Conf_{m}(C^{\circ}_{g,1,{L_{v}}}/\cM_{g,1,{L_{v}}}),\ov{\eta_{C}})\ar[r]\ar[d]&
\pi_{1,\et}(\cM_{g,1,{L_{v}}},\ov{\eta_{M}})\ar[r]\ar[d]\ar@{-->}@/_2pc/[l]_-{s_{P,L_{v}}}&1
\\
1\ar[r]&
\pi_{1,\et}(\Conf_{m}(C_{\ov{\eta}}^{\circ}/\ov{\eta_{M}}),\ov{\eta_{C}})^{\wedge}_{p}\ar[r]&
\pi_{1,\et}'(\Conf_{m}(C^{\circ}_{g,1,\cO_{L_{v}}}/\cM_{g,1,\cO_{L_{v}}}),\ov{\eta_{C}})\ar[r]&
\pi_{1,\et}(\cM_{g,1,\cO_{L_{v}}},\ov{\eta_{M}})\ar[r]\ar@{-->}@/^2pc/[l]^-{s_{P,\cO_{L_{v}}}}&1.
}
\]Here, the notation $\pi_{1,\et}'$ is that of \cite[Expose XIII, \S4]{SGA1}, and the homotopy exact sequence exists on the integral level as $\cO_{L_{v}}$ has residue characteristic $\ne p$.

Therefore, the arithmetic Heisenberg local system $\rho_{Q,p}^{\Heis}$ restricted to $\AMod_{g,1,L_{v}}$ factors through
\[\AMod_{g,1,L_{v}}\thrar\pi_{1,\et}(\cM_{g,1,\cO_{L_{v}}},\ov{\eta_{M}})\rar\Aut_{\cont}((\bB_{m}(\Sigma_{g,1})^{\metab})^{\wedge}_{p}).\]
As $x$ has good reduction at $v$, $x$ extends to a point $\wt{x}\in\cM_{g,1,\cO_{L_{v}}}(\cO_{L_{v}})$, and the rational points give rise to  sections $\Gal(\ov{L_{v}}/L_{v})\rar\pi_{1,\et}(\cM_{g,1,L_{v}},\ov{\eta_{M}})$ and $\Gal(\ov{\cO_{L_{v}}}/\cO_{L_{v}})\rar\pi_{1,\et}(\cM_{g,1,\cO_{L_{v}}},\ov{\eta_{M}})$ that are compatible with each other. Therefore, the restriction of the arithmetic Heisenberg local system at $x$, $\rho_{Q,p}^{\Heis}\rvert_{x}$, restricted to $\Gal(\ov{L_{v}}/L_{v})$ factors through $\Gal(\ov{\cO_{L_{v}}}/\cO_{L_{v}})$, which implies that $c_{Q,x}^{\Heis}$ is unramified at $v$.
\end{proof}
\begin{rmk2}\label{rem:crystalline}
Although we are currently unable to prove it, we also expect that $c_{Q,x}^{\Heis}$ is crystalline under the same hypothesis as Proposition \ref{prop:HeisenbergExtension} with instead $v| p$. This will follow from a generalization of \cite{AIK} to the relative setting, involving crystalline local systems.
\end{rmk2}
\bibliographystyle{alpha}
\bibliography{QuantumLocSys_GQ}

\begin{thebibliography}{DdSMS99}

\bibitem[AI88]{AndersonIhara}
Greg Anderson and Yasutaka Ihara.
\newblock Pro-{$l$} branched coverings of {${\bf P}^1$} and higher circular
  {$l$}-units.
\newblock {\em Ann. of Math. (2)}, 128(2):271--293, 1988.

\bibitem[AIK15]{AIK}
Fabrizio Andreatta, Adrian Iovita, and Minhyong Kim.
\newblock A {$p$}-adic nonabelian criterion for good reduction of curves.
\newblock {\em Duke Math. J.}, 164(13):2597--2642, 2015.

\bibitem[And06]{Andersen}
J\o rgen~Ellegaard Andersen.
\newblock Asymptotic faithfulness of the quantum {${\rm SU}(n)$}
  representations of the mapping class groups.
\newblock {\em Ann. of Math. (2)}, 163(1):347--368, 2006.

\bibitem[BBM19]{BBM}
P.~Belkale, P.~Brosnan, and S.~Mukhopadhyay.
\newblock Hyperplane arrangements and tensor product invariants.
\newblock {\em Michigan Math. J.}, 68(4):801--829, 2019.

\bibitem[BFM23]{BFM}
Prakash Belkale, Najmuddin Fakhruddin, and Swarnava Mukhopadhyay.
\newblock Motivic factorisation of {K}{Z} local systems and deformations of
  representation and fusion rings, 2023.

\bibitem[BFS98]{BFS}
Roman Bezrukavnikov, Michael Finkelberg, and Vadim Schechtman.
\newblock {\em Factorizable sheaves and quantum groups}, volume 1691 of {\em
  Lecture Notes in Mathematics}.
\newblock Springer-Verlag, Berlin, 1998.

\bibitem[BG07]{BG}
Paolo Bellingeri and Eddy Godelle.
\newblock Positive presentations of surface braid groups.
\newblock {\em J. Knot Theory Ramifications}, 16(9):1219--1233, 2007.

\bibitem[BGG17]{BGG}
Paolo Bellingeri, Eddy Godelle, and John Guaschi.
\newblock Abelian and metabelian quotient groups of surface braid groups.
\newblock {\em Glasg. Math. J.}, 59(1):119--142, 2017.

\bibitem[BPS22]{BPS}
Christian Blanchet, Martin Palmer, and Awash Shaukat.
\newblock Heisenberg homology on surface configurations, 2022.

\bibitem[DdSMS99]{Segal}
J.~D. Dixon, M.~P.~F. du~Sautoy, A.~Mann, and D.~Segal.
\newblock {\em Analytic pro-{$p$} groups}, volume~61 of {\em Cambridge Studies
  in Advanced Mathematics}.
\newblock Cambridge University Press, Cambridge, second edition, 1999.

\bibitem[Del89]{Deligne}
Pierre Deligne.
\newblock Le groupe fondamental de la droite projective moins trois points.
\newblock In {\em Galois groups over {${\bf Q}$} ({B}erkeley, {CA}, 1987)},
  volume~16 of {\em Math. Sci. Res. Inst. Publ.}, pages 79--297. Springer, New
  York, 1989.

\bibitem[Fal98]{Faltings}
Gerd Faltings.
\newblock Curves and their fundamental groups (following {G}rothendieck,
  {T}amagawa and {M}ochizuki).
\newblock Number 252, pages Exp. No. 840, 4, 131--150. 1998.
\newblock S\'eminaire Bourbaki. Vol.\ 1997/98.

\bibitem[FWW02]{FWW}
Michael~H. Freedman, Kevin Walker, and Zhenghan Wang.
\newblock Quantum {$\rm SU(2)$} faithfully detects mapping class groups modulo
  center.
\newblock {\em Geom. Topol.}, 6:523--539, 2002.

\bibitem[Gil04]{Gilmer}
Patrick~M. Gilmer.
\newblock Integrality for {TQFT}s.
\newblock {\em Duke Math. J.}, 125(2):389--413, 2004.

\bibitem[GM17]{GM}
Patrick~M. Gilmer and Gregor Masbaum.
\newblock An application of {TQFT} to modular representation theory.
\newblock {\em Invent. Math.}, 210(2):501--530, 2017.

\bibitem[Gro03]{SGA1}
Alexander Grothendieck.
\newblock {\em Rev\^etements \'etales et groupe fondamental ({SGA} 1)},
  volume~3 of {\em Documents Math\'ematiques (Paris) [Mathematical Documents
  (Paris)]}.
\newblock Soci\'et\'e{} Math\'ematique de France, Paris, 2003.
\newblock S\'eminaire de g\'eom\'etrie alg\'ebrique du Bois Marie 1960--61.
  [Algebraic Geometry Seminar of Bois Marie 1960-61], Directed by A.
  Grothendieck, With two papers by M. Raynaud, Updated and annotated reprint of
  the 1971 original [Lecture Notes in Math., 224, Springer, Berlin; MR0354651
  (50 \#7129)].

\bibitem[Iha90]{Ihara}
Yasutaka Ihara.
\newblock {\em Braids, {G}alois groups and some arithmetic functions}.
\newblock ICM-90. Mathematical Society of Japan, Tokyo; distributed outside
  Asia by the American Mathematical Society, Providence, RI, 1990.
\newblock A plenary address presented at the International Congress of
  Mathematicians held in Kyoto, August 1990.

\bibitem[Lan24]{Landesman}
Aaron Landesman.
\newblock Invariance of the tame fundamental group under base change between
  algebraically closed fields.
\newblock {\em Essent. Number Theory}, 3(1):1--18, 2024.

\bibitem[Loo12]{Looijenga}
Eduard Looijenga.
\newblock The {KZ} system via polydifferentials.
\newblock In {\em Arrangements of hyperplanes---{S}apporo 2009}, volume~62 of
  {\em Adv. Stud. Pure Math.}, pages 189--231. Math. Soc. Japan, Tokyo, 2012.

\bibitem[Lys22]{Lysenko}
Sergey Lysenko.
\newblock Towards canonical representations of finite {H}eisenberg groups.
\newblock {\em Bull. Soc. Math. France}, 150(3):569--577, 2022.

\bibitem[Mat96]{Matsumoto}
Makoto Matsumoto.
\newblock Galois representations on profinite braid groups on curves.
\newblock {\em J. Reine Angew. Math.}, 474:169--219, 1996.

\bibitem[MN22]{MinamideNakamura}
Arata Minamide and Hiroaki Nakamura.
\newblock The automorphism groups of the profinite braid groups.
\newblock {\em Amer. J. Math.}, 144(5):1159--1176, 2022.

\bibitem[Moc99]{Mochizuki}
Shinichi Mochizuki.
\newblock The local pro-{$p$} anabelian geometry of curves.
\newblock {\em Invent. Math.}, 138(2):319--423, 1999.

\bibitem[NS07]{NikolovSegal}
Nikolay Nikolov and Dan Segal.
\newblock On finitely generated profinite groups. {I}. {S}trong completeness
  and uniform bounds.
\newblock {\em Ann. of Math. (2)}, 165(1):171--238, 2007.

\bibitem[Oda97]{Oda}
Takayuki Oda.
\newblock Etale homotopy type of the moduli spaces of algebraic curves.
\newblock In {\em Geometric {G}alois actions, 1}, volume 242 of {\em London
  Math. Soc. Lecture Note Ser.}, pages 85--95. Cambridge Univ. Press,
  Cambridge, 1997.

\bibitem[RM23]{DeRenziMartel}
Marco~De Renzi and Jules Martel.
\newblock Homological {C}onstruction of {Q}uantum {R}epresentations of
  {M}apping {C}lass {G}roups, 2023.

\bibitem[{Sta}18]{stacks}
The {Stacks Project Authors}.
\newblock \textit{Stacks Project}.
\newblock \url{https://stacks.math.columbia.edu}, 2018.

\end{thebibliography}
\end{document}